\documentclass[leqno,english]{amsproc}
\usepackage{hyperref}
\usepackage{amsthm}
\usepackage{amssymb}
\usepackage[square,numbers]{natbib}
\usepackage[all]{xy}

\numberwithin{subsubsection}{section}

\theoremstyle{plain}

\newtheorem*{mainthm}{Theorem}

\newtheorem{thm}[subsubsection]{Theorem}
\newtheorem{lemm}[subsubsection]{Lemma}
\newtheorem{prop}[subsubsection]{Proposition}

\DeclareMathOperator{\A}{\mathcal{A}}
\DeclareMathOperator{\C}{\mathcal{C}}

\DeclareMathOperator{\Op}{\mathcal{O}}

\DeclareMathOperator{\Mor}{\mathrm{Mor}}
\DeclareMathOperator{\Hom}{\mathrm{Hom}}
\DeclareMathOperator{\End}{\mathit{End}}

\DeclareMathOperator{\Indec}{\mathrm{Indec}}

\DeclareMathOperator*{\colim}{\mathrm{colim}}
\DeclareMathOperator*{\im}{\mathrm{Im}}

\DeclareMathOperator{\NN}{\mathbb{N}}
\DeclareMathOperator{\ZZ}{\mathbb{Z}}
\DeclareMathOperator{\FF}{\mathbb{F}}

\DeclareMathOperator{\kk}{\Bbbk}

\DeclareMathOperator{\AOp}{\mathsf{A}}

\DeclareMathOperator{\COp}{\mathsf{C}}
\DeclareMathOperator{\DOp}{\mathsf{D}}
\DeclareMathOperator{\EOp}{\mathsf{E}}

\DeclareMathOperator{\GOp}{\mathsf{G}}

\DeclareMathOperator{\IOp}{\mathsf{I}}

\DeclareMathOperator{\LOp}{\mathsf{L}}
\DeclareMathOperator{\POp}{\mathsf{P}}
\DeclareMathOperator{\QOp}{\mathsf{Q}}

\DeclareMathOperator{\ad}{\mathit{ad}}
\DeclareMathOperator{\id}{\mathit{id}}

\title[Koszul duality complexes for iterated loop space cohomology]{Koszul duality complexes for the cohomology of iterated loop spaces of spheres}

\author{Benoit Fresse}
\date{6 January 2013}
\address{UMR 8524 de l'Universit\'e de Lille 1 - Sciences et Technologies - et du CNRS\\
Cit\'e Scientifique -- B\^atiment M2\\
F-59655 Villeneuve d'Ascq C\'edex (France)}
\email{Benoit.Fresse@math.univ-lille1.fr}
\urladdr{http://math.univ-lille1.fr/\~{ }fresse}
\thanks{Research supported in part by grant ANR-06-JCJC-0042 ``OBTH''}

\begin{document}

\begin{abstract}
The goal of this article is to make explicit a structured complex whose homology computes the cohomology of the $p$-profinite completion
of the $n$-fold loop space of a sphere
of dimension $d = n-m<n$.
This complex is defined purely algebraically, in terms of characteristic structures
of $E_n$-operads.
Our construction involves: the free complete algebra in one variable associated to any $E_n$-operad;
and an element in this free complete algebra, which is associated to a morphism
from the operad of $L_\infty$-algebras
to an operadic suspension of our $E_n$-operad.
We deduce our main theorem from: a connection between the cohomology of iterated loop spaces
and the cohomology of algebras over $E_n$-operads;
and a Koszul duality result for $E_n$-operads.
\end{abstract}

\maketitle

\section*{Introduction}

The goal of this article is to make explicit a structured complex whose homology computes the cohomology of the $p$-profinite completion
of the $n$-fold loop space of a sphere of dimension $d = n-m<n$.

The homology of the $n$-fold loop space on the $n$-sphere $\Omega^n S^n$ was first determined by Araki-Kudo in~\cite{ArakiKudo},
in the $\FF_2$-coefficient case,
by using a construction of natural homological operations,
later called the Araki-Kudo-Dyer-Lashof operations,
associated to iterated loop spaces.
The definition of these operations was generalized by F. Cohen in~\cite{Cohen}.
The outcome of Cohen's construction is a computation of the homology of the $n$-fold space of a suspension $\Omega^n\Sigma^n X$
as a functor on the homology of the space $X$.
The work of F. Cohen also gives the connection between the algebra of the Araki-Kudo-Dyer-Lashof operations
and the equivariant homology of the Boardman-Vogt little $n$-cubes operad $\COp_n$.
The homology of the $n$-fold loop space $\Omega^n S^{n-m}$, where $m>0$,
has been computed by Hunter in~\cite{Hunter} in the case $m=1$,
and by Cohen-Peterson in~\cite{CohenPeterson} in the case $n=\infty$ and $m=2$.
Little is known on this homology otherwise.

The goal of this paper is not to give a computation method, but rather to give a chain complex which makes explicit
the connection between the cohomology of such iterated loop space~$\Omega^n S^{n-m}$
and homotopical structures attached to~$E_n$-operads.
For our purpose, we consider $E_n$-operads in the category of differential graded modules (dg-modules for short)
and we use the noun $E_n$-algebra to refer to any category of algebras in dg-modules
associated to an $E_n$-operad. In our setting, an $E_n$-operad is precisely defined as an operad which is weakly equivalent
to the chain operad of Boardman-Vogt' little $n$-cubes $\COp_n$.
Let us explain the result of our construction.

\medskip
Most models of $E_n$-operads
come in nested sequences
\begin{equation}\label{EnOperadNestedSequence}
\EOp_1\subset\cdots\subset\EOp_{n-1}\subset\EOp_n\subset\cdots
\end{equation}
such that $\EOp_{\infty} = \colim_n\EOp_n$ is an $E_{\infty}$-operad, an operad weakly equivalent to the operad of associative and commutative algebras $\COp$.
The $E_1$-operads, which form the initial terms of nested sequences of this form (\ref{EnOperadNestedSequence}),
are weakly equivalent to the operad
of associative non-commutative algebras $\AOp$.
Intuitively, an $E_1$-algebra structure is interpreted as a structure which is homotopy associative in the strong sense but non-commutative,
while an $E_\infty$-algebra is fully homotopy (associative and) commutative.
The $E_n$-algebra structures, which are intermediate between $E_1$ and $E_\infty$,
are interpreted as structures which are homotopy associative in the strong sense
and homotopy commutative up to some degree. A commutative algebra naturally inherits an $E_n$-algebra structure,
for all $n = 1,2,\dots,\infty$.

We also consider the operad $\LOp_{\infty}$ associated to the category of $L_\infty$-algebras, where an $L_{\infty}$-algebra consists of a dg-module $L$
equipped with an $r$-ary Lie bracket $\lambda_r: L^{\otimes r}\rightarrow L$, for every $r\geq 2$,
so that higher analogues of the usual antisymmetry and Jacobi relations hold in $L$ (see~\cite{HinichSchechtmanLie}).
We have an operad morphism $\phi_m^{\sharp}: \Lambda\LOp_{\infty}\rightarrow\Lambda^m\EOp_m$,
for each $m\in\NN$,
where $\Lambda$ refers to an operadic suspension operation (see our recollections in~\S\ref{Background:Suspensions}).
We prove the existence and homotopy uniqueness of such morphisms in~\cite{EnKoszulDuality,OperadMaps}
by using obstruction theory methods for operads (we give an outline of this construction in~\S\ref{KoszulDualityReview}).
We refer to~\cite{KontsevichMotives} for another construction of these morphisms, relying on Kontsevich's complex of semi-algebraic forms
and on a particular model of $E_n$-operad, the Fulton-MacPherson operad,
introduced by Getzler-Jones in~\cite{GetzlerJones}.

The existence of this morphism $\phi_m^{\sharp}: \Lambda\LOp_{\infty}\rightarrow\Lambda^m\EOp_m$
reflects the existence of a restriction functor from $E_m$-algebras to $L_\infty$-algebras
and has been used in the literature in the case $m=2$,
in combination with the Deligne conjecture,
in order to give a new proof of Kontsevich's formality theorem
in deformation quantization (see~\cite{KontsevichMotives}).
Recall that the homology of the little $m$-discs operad is identified with the Gerstenhaber operad $\GOp_m$,
which is a composite of the commutative operad $\COp$,
and of the $m-1$-fold operadic desuspension of the Lie operad $\Lambda^{1-m}\LOp$.
The morphism $\phi_m^{\sharp}: \Lambda\LOp_{\infty}\rightarrow\Lambda^m\EOp_m$ corresponds, at the homology level,
to the canonical embedding morphism from the (suspension of the) operad of Lie algebras $\LOp$
towards the ($m$-fold suspension of the) Gerstenhaber operad $\GOp_m$.
In~\cite{OperadMaps}, we precisely prove that the morphism $\phi_m^{\sharp}$ is uniquely characterized by this induced action in homology,
up to a right homotopy in the category of operads,
and as such, is characteristic of the structure of an~$E_m$-operad.

We form, for any $n\geq m$, the free complete $\EOp_n$-algebra $\widehat{\EOp}_n(x)$
on one generator $x$ of degree~$-m$. The elements of this algebra~$\widehat{\EOp}_m(x)$
are just formal power series $\sum_{r=1}^{\infty} \xi_r(x,\dots,x)$
where $\xi_r$ is an $r$-ary operation of the operad $\EOp_m$.
We consider, after~\cite{ChapotonLivernetHopf}, a composition operation $\circ: \widehat{\EOp}_n(x)\otimes\widehat{\EOp}_n(x)\rightarrow\widehat{\EOp}_n(x)$,
defined termwise by sums of operadic substitution operations
\begin{equation*}
\xi_r(x,\dots,x)\circ\zeta_s(x,\dots,x) = \sum_{i=1}^{r} \xi_r(x,\dots,\zeta_s(x,\dots,x),\dots,x)
\end{equation*}
running over all factors $x$ of the operadic monomial $\xi_r(x,\dots,x)$.
We check in~\S\ref{FreeCompleteAlgebras:LinfinityMorphisms} that the existence of a morphism $\phi_m^{\sharp}: \Lambda\LOp_{\infty}\rightarrow\Lambda^m\EOp_m$
is formally equivalent to the existence of an element $\omega_m\in\widehat{\EOp}_m(x)$,
of degree $-1-m$,
and satisfying the relation $\delta(\omega_m) = \omega_m\circ\omega_m$ in $\widehat{\EOp}_m(x)$,
where $\delta$ refers to the canonical differential of the dg-module $\widehat{\EOp}_m(x)$.
In~\cite{EnKoszulDuality}, we actually give the construction of such an element $\omega_m$,
by using obstruction theory methods,
rather than the construction of an equivalent operad morphism $\phi_m^{\sharp}: \Lambda\LOp_{\infty}\rightarrow\Lambda^m\EOp_m$.
(We go back to this construction in~\S\ref{SphereKoszulDuality:CobarAugmentations}.)
We embed this element $\omega_m$ into $\widehat{\EOp}_n(x)$, for any $n\geq m$.

Our result reads:

\begin{mainthm} 
We take a finite primary field $\FF = \FF_p$ of characteristic $p>0$
as ground ring.
We consider the twisted complex $(\widehat{\EOp}_n(x),\partial_m)$ defined by adding the map $\partial_m(\xi) = \xi\circ\omega_m$
to the natural differential of the free complete $\EOp_n$-algebra $\widehat{\EOp}_n(x)$.

We have an identity
\begin{equation*}
H^*(\widehat{\Omega^n S^{n-m}}) = H_*(\widehat{\EOp}_n(x)^{\vee},\partial_m^{\vee})
\end{equation*}
between the continuous cohomology with $\FF_p$ coefficients of the $p$-profinite completion of $\Omega^n S^{n-m}$
and the homology of the dual complex of $(\widehat{\EOp}_n(x),\partial_m)$
in $\FF_p$-modules (we take the continuous dual with respect to the topology of power series).
\end{mainthm}

The profinite completion of the theorem 
is defined by the limit of a diagram $\Omega^n S^{n-m}_\alpha$
where $S^{n-m}_\alpha$ ranges over homotopy $p$-finite approximations of $S^{n-m}$ (see~\cite{Morel}).
This formulation is slightly abusive because the $p$-profinite completion
of $\Omega^n S^{n-m}$ would be the limit of a diagram of homotopy $p$-finite approximations of the space $\Omega^n S^{n-m}$
itself. Nevertheless results of~\cite{Shipley}
imply that these limits coincide.

\medskip
Let us explain the background of the proof of our theorem.
Throughout this paper, we use the noun space for simplicial set. To any simplicial set $X$, we associate the normalized complex $N^*(X)$
with coefficients in a ground ring $\FF$
of our base category of dg-modules. In what follows, we assume that this ground ring is a finite primary field $\FF = \FF_p$
of characteristic $p>0$.
Since we deal with pointed spaces, we use the reduced complex $\bar{N}^*(X) = \ker(N^*(X)\rightarrow\FF)$
in our constructions rather than the whole complex $N^*(X)$.
To handle structures without unit (like these reduced complexes $\bar{N}^*(X)$),
we also tacitely assume that we deal with non-unitary operads (in the sense of~\cite{OperadModules}),
with no operation in arity zero.

The cochain algebras $\bar{N}^*(X)$ inherit an $E_\infty$-algebra structure. This result is establish in~\cite{HinichSchechtman}
by an abstract argument. We also have an explicit construction of this $E_\infty$-algebra structure, given in~\cite{BergerFresse,McClureSmith},
and involving combinatorial models of~$E_\infty$-operads
equipped with a filtration of the form~(\ref{EnOperadNestedSequence}).
In this situation, we can use an obvious restriction process
to make explicit an~$E_n$-algebra structure on the cochain complex~$\bar{N}^*(X)$,
for every $n\in\NN$.

We have a natural homology theory $H_*^{\POp}(-)$ associated to any operad $\POp$.
We establish in~\cite{IteratedBar} that the $\EOp_n$-homology $H_*^{\EOp_n}(\bar{N}^*(X))$ of a cochain algebra~$A = \bar{N}^*(X)$
determines the (continuous) cohomology of the (homotopy $p$-profinite completion) of the $n$-fold loop space on $X$.
We prove in~\cite{EnKoszulDuality} that the $E_n$-operads are Koszul, for all $n = 1,2,\dots$, and this result implies that the homology of an $E_n$-algebra
is determined by an explicit chain complex, given as a quasi-cofree coalgebra over the Koszul dual cooperad
of any model of~$E_n$-operad~$\EOp_n$.

The Koszul dual of an $E_n$-operad $\EOp_n$ is actually identified with the cooperad such that $\DOp_n = \Lambda^{-n}\EOp_n^{\vee}$,
where $\EOp_n^{\vee}$ denotes the dual cooperad of the operad~$\EOp_n$
in the category of dg-modules. (Recall that $\Lambda$ denotes an operadic suspension operation.)
We deduce from this statement that the Koszul duality complex $C^{\EOp_n}_*(A)$
which determines the homology $H_*^{\EOp_n}(A)$
of an $\EOp_n$-algebra $A$
has the form $C^{\EOp_n}_*(A) = (\DOp_n(A),\partial)$,
where
$\DOp_n(A)$ is the connected $\DOp_n$-cofree coalgebra on~$A$
and
$\partial$ is a twisting coderivation determined by the $E_n$-algebra structure of~$A$ (we review the details of this construction
in~\S\ref{KoszulDualityReview}).

The main purpose of this article is to determine the twisting coderivation $\partial$
associated to this Koszul duality complex in the case of the cochain algebra of a sphere $A = \bar{N}^*(S^{n-m})$,
and to check that the obtained complex $C^{\EOp_n}_*(\bar{N}^*(S^{n-m}))$
is identified with the continuous dual of the dg-module $(\widehat{\EOp}_n(x),\partial_m)$
considered in our theorem.

\medskip
In~\S\ref{KoszulDualityReview}, we briefly review the definition of the homology $H^{\POp}_*(A)$ associated to an operad $\POp$,
the applications of the Koszul duality of operads to the definition of complexes $C^{\POp}_*(A)$
computing $H^{\POp}_*(A)$,
and we explain how to apply this construction and the result of~\cite{EnKoszulDuality}
to $E_n$-operads.
In~\S\ref{FreeCompleteAlgebras}, we explain the definition of the free complete algebras $\widehat{\EOp}_n(x)$
and we check that this free complete algebra is dual to a cofree coalgebra $\DOp(C)$
on a dg-module of rank one $C$.
In~\S\ref{SphereKoszulDuality}, we study the applications of operadic Koszul duality to the cochain algebra of spheres $\bar{N}^*(S^{n-m})$
and we determine the form of the complex $C^{\EOp_n}_*(\bar{N}^*(S^{n-m}))$
from which we obtain the result of our theorem.

In the constructions of these sections~\S\S\ref{KoszulDualityReview}-\ref{SphereKoszulDuality}, we deal with a particular instance of $E_n$-operad,
already used in the proof of~\cite{EnKoszulDuality}, and which arises from a certain filtration of an $E_\infty$-operad
introduced by Barratt-Eccles in~\cite{BarrattEccles}.
To complete our results, we establish in~\S\ref{HomotopyInvariance} that the homotopy type of the complex $(\widehat{\EOp}_n(x)^{\vee},\partial_m^{\vee})$
does not depend on the choice of the $E_n$-operad $\EOp_n$ and of the element $\omega_m\in\widehat{\EOp}_n(x)$
associated to a morphism $\phi_m^{\sharp}: \Lambda\LOp_{\infty}\rightarrow\Lambda^n\EOp_n$.
This verification proves our claim that our chain complex, computing the cohomology of iterated loop spaces of spheres,
is determined by homotopy structures
attached to $E_n$-operads.

Before starting our main matter, we devote a preliminary section to a comprehensive review of the background of our constructions.

\section*{Background and general conventions}\label{Background}

The main purpose of this preliminary is to explain our conventions. We refer to~\cite{CylinderOperads}
for a more detailed account of this background,
with the same conventions as in the present article, and further bibliographical references.
To begin with, we specify which base category of dg-modules
is considered throughout the article.

\subsubsection{The category of dg-modules}\label{Background:DGModules}
In~\cite{CylinderOperads}, we work in a category of dg-modules over an arbitrary ground ring $\kk$.
In this paper, we will assume that the ground ring is a field $\kk = \FF$.
This assumption simplifies the statement of some results on dg-modules.
For the moment, we do not make any requirement on the characteristic of this ground field, but when we tackle topological applications,
we have to assume that this characteristic is non-zero.

For us a dg-module refers to a $\ZZ$-graded $\FF$-module $C$ equipped with an internal differential $\delta: C\rightarrow C$
that decreases degrees by one. We use the letter $\C$ to refer to this category of dg-modules
which we take as base category for all our constructions.

The category of dg-modules is equipped with its standard tensor product
whose symmetry isomorphism $\tau: C\otimes D\rightarrow D\otimes C$
involves a sign.
This sign is determined by the standard sign convention of differential graded algebra.
The notation $\pm$ is used to represent any sign arising from an application of this convention.
The ground ring $\FF$, viewed as a dg-module of rank~$1$ concentrated in degree~$0$,
defines the unit object for the tensor product of dg-modules.

The morphism sets of any category $\A$ are denoted by $\Mor_{\A}(A,B)$.
The morphisms of dg-modules are the degree and differential preserving morphisms of $\FF$-modules $f: C\rightarrow D$.
The category of dg-modules comes also equipped with internal hom-objects $\Hom_{\C}(C,D)$
characterized by the adjunction relation
$\Mor_{\C}(K\otimes C,D) = \Mor_{\C}(K,\Hom_{\C}(C,D))$
with respect to the tensor product $\otimes: \C\times\C\rightarrow\C$.
Recall briefly that a homogeneous element of $\Hom_{\C}(C,D)$
is just a morphism of $\FF$-modules $f: C\rightarrow D$
raising degrees by $d = \deg(f)$.
The differential of a homomorphism $f\in\Hom_{\C}(C,D)$
is defined by the (graded) commutator $\delta(f) = \delta\cdot f - \pm f\cdot\delta$,
where we consider the internal differentials of $C$ and $D$.
The elements of the dg-hom $\Hom_{\C}(C,D)$
are called homomorphisms to be distinguished from the actual morphisms of dg-modules,
which represent the elements of the morphism set $\Mor_{\C}(C,D)$.

The dual of a dg-module $C$ is the dg-module such that $C^{\vee} = \Hom_{\C}(C,\FF)$,
where we again view the ground ring $\FF$
as a dg-module of rank~$1$ concentrated in degree~$0$.

\subsubsection{Operads and cooperads}\label{Background:OperadsCooperads}
We always consider operads (and cooperads) in the category of dg-modules.
We adopt the notation $\Op$ for the category of operads.
We use the notation $\IOp$ for the initial object of this category, the operad such that $\IOp(1) = \FF$
and $\IOp(r) = 0$ for $r\not=1$.

We assume throughout this article that an operad $\POp$ satisfies $\POp(0) = 0$. (In our work, we often use the expression of non-unitary operad
to refer to this requirement.)
We assume in some cases that our operad also satisfies $\POp(1) = \FF$ in addition to $\POp(0) = 0$.
We say in this situation that $\POp$ is connected as an operad.
Cooperads $\DOp$ are always assumed to satisfy the connectedness assumptions $\DOp(0) = 0$ and $\DOp(1) = \FF$
in order to avoid difficulties with infinite sums in coproducts.
We use the notation $\Op_0$ (respectively $\Op_1$) for the category of non-unitary (respectively connected) operads in dg-modules.
We call $\Sigma_*$-object the structure, underlying an operad, formed by a collection $M = \{M(r)\}_{r\in\NN}$
where $M(r)$ is a dg-modules equipped with an action of the symmetric group on $r$ letters $\Sigma_r$.

A connected operad inherits a canonical augmentation $\epsilon: \POp\rightarrow\IOp$
which is the identity in arity $r=1$
and is trivial otherwise.
We adopt the notation $\widetilde{\POp}$ for the augmentation ideal of any operad $\POp$ equipped with an augmentation $\epsilon: \POp\rightarrow\IOp$.
We have in the connected case $\widetilde{\POp}(0) = \widetilde{\POp}(1) = 0$ and $\widetilde{\POp}(r) = \POp(r)$ for all $r\geq 2$.
Any connected cooperad $\DOp$ has a coaugmentation $\eta: \IOp\rightarrow\DOp$
and a coaugmentation coideal $\widetilde{\DOp}$
defined like the augmentation and the augmentation ideal
of a connected operad.

\subsubsection{Suspensions}\label{Background:Suspensions}
For any $d\in\ZZ$, we adopt the notation $\FF[d]$ for the free graded $\FF$-module of rank $1$
concentrated in degree $d$.
The suspension of a dg-module $C$ is the dg-module such that $\Sigma C = \FF[1]\otimes C$

The operadic suspension of an operad $\POp$,
is an operad $\Lambda\POp$
characterized by the commutation relation $\Lambda\POp(\Sigma C) = \Sigma\POp(C)$
at the level of free algebras.
Basically,
the operad $\Lambda\POp$ is defined arity-wise by the tensor products $\Lambda\POp(n) = \FF[1-n]\otimes\POp(n)^{\pm}$
where the notation $\pm$
refers to a twist of the natural $\Sigma_n$-action on $\POp(n)$
by the signature of permutations.
We have an operadic suspension operation defined in the same way on cooperads.

We may also apply the suspension of dg-modules arity-wise to any operad $\POp$ (or cooperad)
in order to produce a $\Sigma_*$-object $\Sigma\POp$
such that $\Sigma\POp(r) = \FF[1]\otimes\POp(r)$,
but this suspended $\Sigma_*$-object does not inherit an operad structure.

\subsubsection{Algebra and coalgebra categories}\label{Background:AlgebrasCoalgebras}
The category of algebras associated to an operad $\POp$
is denoted by~${}_{\POp}\C$.
The free $\POp$-algebra associated to a dg-module $C$ is denoted by $\POp(C)$.
Recall simply that this free $\POp$-algebra is defined by the dg-module of generalized symmetric tensors
$\POp(C) = \bigoplus_{r=0}^{\infty} (\POp(r)\otimes C^{\otimes r})_{\Sigma_r}$,
where the $\Sigma_r$-quotient identifies tensor permutations
with the action of permutations on~$\POp(r)$.
Note that we can remove the $0$ term from this expansion since we assume $\POp(0) = 0$.
We use the notation $p(x_1,\dots,x_r)$, where $p\in\POp(r)$ and $x_1\otimes\dots\otimes x_r\in A^{\otimes r}$,
for the element of $\POp(C)$
defined by the tensor $p\otimes(x_1\otimes\dots\otimes x_r)$.
We may also use the graphical representation
\begin{equation}
\vcenter{\xymatrix@H=6pt@W=3pt@M=2pt@!R=1pt@!C=1pt{ x_1\ar[dr]\ar@{.}[r] & \cdots\ar@{.}[r]\ar@{}[d]|{\displaystyle{\cdots}} & x_r\ar[dl] \\
\ar@{.}[r] & *+<8pt>[F]{p}\ar[d]\ar@{.}[r] & \\
\ar@{.}[r] & 0\ar@{.}[r] & }}
\end{equation}
for this element $p(x_1,\dots,x_r)\in\POp(C)$.

In the case of a cooperad $\DOp$, we use the notation $\DOp(C)$ for the connected cofree $\DOp$-coalgebra
such that $\DOp(C) = \bigoplus_{r=0}^{\infty} (\DOp(r)\otimes C^{\otimes r})_{\Sigma_r}$.
In this construction, we take the same generalized symmetric algebra functor as in the operad case.
Therefore, we can also adopt the functional and graphical representation of elements of free algebras of operads to denote the elements of~$\DOp(C)$.
Note that $\DOp(C)$ does not agree with the standard cofree connected coalgebra construction because we take coinvariants instead of invariants.
But we generally apply the definition of this object $\DOp(C)$ to cooperads $\DOp$
equipped with a free $\Sigma_*$-module structure, and as a consequence,
the result of our construction
does not change if we replace coinvariants by invariants.

The coaugmentation morphism of a cooperad yields an embedding $\eta: C\hookrightarrow\DOp(C)$
which identifies the dg-module $C$ with a summand of~$\DOp(C)$.
Similarly, if $\POp$ is an augmented operad, then we have a split embedding $\eta: C\hookrightarrow\POp(C)$,
yielded by the unit of the operad~$\POp$.
The morphism $\epsilon: \POp(C)\rightarrow C$ induced by the augmentation of~$\POp$
gives the retraction such that $\epsilon\eta = \id$.

\subsubsection{Model structures}\label{Background:ModelCategories}
The category of dg-modules $\C$ is equipped with its standard model structure in which the weak-equivalences are the morphisms
which induce an isomorphism in homology, the fibrations are the degreewise surjections (see~\cite[\S 2.3]{Hovey}).
The assumption that the ground ring is a field implies that the class of cofibrations of the category of dg-modules
is identified with the class of injective morphisms
of dg-modules. In particular, any dg-module forms a cofibrant object in our setting.

The category of non-unitary operads $\Op_0$ inherits a model structure with as weak-equivalences (respectively, fibrations) the morphisms $f: \POp\rightarrow\QOp$
which form a weak-equivalence (respectively, fibration) of dg-modules $f: \POp(r)\rightarrow\QOp(r)$
in each arity $r\in\NN$~(see~\cite{HinichHomotopy} or~\cite{BergerMoerdijk}).
The cofibrations of the category of operads are characterized by the right-lifting-property with respect to the acyclic fibrations.

The category of $\Sigma_*$-objects in dg-modules also inherits a model structure (like any category of modules over a ring~\cite[\S 2.3]{Hovey}).
The weak-equivalences (respectively, fibrations) of this category are created arity-wise in the category of dg-modules,
as in the operad case.

Recall that any operad cofibration with a cofibrant operad as domain defines a cofibration of $\Sigma_*$-objects (see again~\cite{BergerMoerdijk}
and the corrections in~\cite{BergerMoerdijkColoured}),
but the converse implication does not hold. Therefore we say that an operad $\POp$ is $\Sigma_*$-cofibrant when its unit morphism $\eta: \IOp\rightarrow\POp$
defines a cofibration in the category of $\Sigma_*$-objects.

The category of algebras over a $\Sigma_*$-cofibrant operad $\POp$ inherits a semi-model structure
in the sense that the axioms of model categories are satisfied for $\POp$-algebras
when we restrict ourselves to acyclic cofibrations with a cofibrant domain
in Quillen's lifting axiom (M4.ii),
respectively to morphisms with a cofibrant domain in Quillen's factorization axioms (M5.i-ii).
This is enough for the usual constructions of homotopical algebra. (We refer to~\cite{OperadModules} and for a comprehensive account on this background.)
The weak-equivalences (respectively, fibrations) of $\POp$-algebras are the morphisms of $\POp$-algebras $f: A\rightarrow B$
which define a weak-equivalence (respectively, a fibration)
in the category of dg-modules. The cofibrations are characterized by the right-lifting-property with respect to acyclic fibrations.

\subsubsection{Twisted objects}\label{Background:TwistedObjects}
In certain constructions, a homomorphism $\partial\in\Hom_{\C}(C,C)$ of degree $-1$
is added to the internal differential of a dg-module $C$
in order to produce a new dg-module, with the same underlying graded $\FF$-module as $C$,
but with the map $\delta+\partial: C\rightarrow C$
as differential.
The relation of differential $(\delta+\partial)^2 = 0$
is equivalent to the equation $\delta(\partial) + \partial^2 = 0$
in $\Hom_{\C}(C,C)$.
In this situation, we say that $\partial$ is a twisting homomorphism and we use the pair $(C,\partial)$ as notation
for the new dg-module defined by the addition of $\partial$
to the internal differential of $C$.

In the case of an algebra $A$ over an operad $\POp$, we assume that the twisting homomorphism
satisfies the derivation relation
\begin{equation*}
\partial(p(a_1,\dots,a_r)) = \sum_{i=1}^{r} p(a_1,\dots,\partial(a_i),\dots,a_r)
\end{equation*}
in order to ensure that the twisted dg-module $(A,\partial)$
inherits a $\POp$-algebra structure
from $A$.
We then say that $\partial: A\rightarrow A$
is a twisting derivation.

We call quasi-free $\POp$-algebra the twisted $\POp$-algebras $(A,\partial)$
associated to a free $\POp$-algebra $A = \POp(C)$.
We have a one-one correspondence between derivations on free $\POp$-algebras $\partial: \POp(C)\rightarrow\POp(C)$
and homomorphisms $\gamma: C\rightarrow\POp(C)$.
We use the notation $\partial_{\gamma}$
for the derivation associated to $\gamma$.
In one direction, we have $\gamma = \partial_{\gamma}|_C$.
In the other direction, the derivation $\partial_{\gamma}: \POp(C)\rightarrow\POp(C)$
is determined from $\gamma: C\rightarrow\POp(C)$
by the derivation formula $\partial_{\gamma}(p(c_1,\dots,c_r)) = \sum_{i=1}^{r} p(c_1,\dots,\gamma(c_i),\dots,c_r)$,
for all $p(c_1,\dots,c_r)\in\POp(C)$.

For a free $\POp$-algebra $\POp(C)$, a derivation $\partial_{\gamma}: \POp(C)\rightarrow\POp(C)$
satisfies the equation of twisting homomorphisms
if and only if the associated homomorphism $\gamma: C\rightarrow\POp(C)$
satisfies the equation $\delta(\gamma) + \partial_{\gamma}\cdot\gamma = 0$
in $\Hom_{\C}(C,\POp(C))$.

We also apply the definition of twisted dg-modules to coalgebras over cooperads.
We have a notion of quasi-cofree $\DOp$-coalgebra,
dual to the notion of quasi-free algebra over an operad,
consisting of a twisted dg-module $(D,\partial)$ such that $D = \DOp(C)$ is a cofree $\DOp$-coalgebra
and $\partial = \partial_\alpha: \DOp(C)\rightarrow\DOp(C)$ is a $\DOp$-coalgebra coderivation
uniquely determined by a homomorphism $\alpha: \DOp(C)\rightarrow C$
such that $\delta(\alpha) + \alpha\cdot\partial_\alpha = 0$
in $\Hom_{\C}(\DOp(C),C)$.
In the context of coalgebras,
we determine the twisting coderivation $\partial_\alpha$ associated to a homomorphism $\alpha$
by a graphical expression of the form:
\begin{multline}
\partial_{\alpha}
\left\{\vcenter{\xymatrix@H=6pt@W=3pt@M=2pt@!R=1pt@!C=1pt{ x_1\ar[dr]\ar@{.}[r] & \cdots\ar@{.}[r]\ar@{}[d]|{\displaystyle{\cdots}} & x_r\ar[dl] \\
\ar@{.}[r] & *+<8pt>[F]{c}\ar[d]\ar@{.}[r] & \\
\ar@{.}[r] & 0\ar@{.}[r] & }}\right\}
= \sum_{i}
\pm\left\{\vcenter{\xymatrix@H=6pt@W=3pt@M=2pt@!R=1pt@!C=1pt{ x_1\ar[drr]\ar@{.}[r] & \cdots\ar@{.}[r] &
\alpha(x_i)\ar[d]\ar@{.}[r] & \cdots\ar@{.}[r] & x_r\ar[dll] \\
\ar@{.}[rr] && *+<8pt>[F]{c}\ar@{.}[rr]\ar[d] && \\
\ar@{.}[rr] && 0\ar@{.}[rr] && }}\right\}
\\
+ \sum_{\substack{\tau\in\Theta_2(r)\\ \rho_{\tau}(c)}}
\pm\left\{\vcenter{\xymatrix@H=6pt@W=4pt@M=2pt@R=8pt@C=4pt{ &&
\save [].[rrd]!C *+<4pt>[F-,]\frm{}*+<6pt>\frm{\{}*+<0pt>\frm{\}} \restore\ar@{}[]!L-<8pt,0pt>;[d]!L-<8pt,0pt>_{\displaystyle\alpha}
x_*\ar@{.}[r]\ar[dr] & \cdots\ar@{.}[r] & x_*\ar[dl] && \\
x_*\ar@{.}[r]\ar@/_6pt/[drrr] & \cdots && *+<6pt>[F]{c_*}\ar@{.}[r]\ar@{.}[l]\ar[d] && \cdots\ar@{.}[r] & x_*\ar@/^6pt/[dlll] \\
\ar@{.}[rrr] &&& *+<6pt>[F]{c_*}\ar[d]\ar@{.}[rrr] &&& \\
\ar@{.}[rrr] &&& 0\ar@{.}[rrr] &&& }}\right\}.
\end{multline}
The notation~$\Theta_2(r)$
refers to the category of trees with two vertices.
The notation~$\rho_{\tau}(c)$ refers to a cooperadic coproduct along any tree~$\tau\in\Theta_2(r)$
of any element $c\in\DOp(r)$.
The expressions $c_*$ inside the sum refer to the factors of this cooperadic coproduct (see~\cite[Proposition 4.1.3]{CylinderOperads} for details).
The second sum ranges over all trees $\tau\in\Theta_2(r)$
and, for each $\tau\in\Theta_2(r)$, over all terms of the expansion of $\rho_{\tau}(c)$.
In the first sum, we just consider the restriction of $\alpha$
to the summand $C\subset\DOp(C)$
and we sum over all inputs of $c\in\DOp(r)$.

\section{Koszul duality and operadic homology}\label{KoszulDualityReview}

The purpose of this section is to review applications of the bar duality of operads for the definition of homology theories
associated to categories algebras over an operad.
The Koszul duality refers to a particular case of operadic bar duality, where a cobar model of an operad $\POp$
can be determined from a presentation of the operad $\POp$. The $E_n$-operads are not Koszul in this sense,
but the homology of an $E_n$-operad is identified with the $n$-Gerstenhaber operad,
which forms an instance of a self-dual Koszul operad~\cite{GetzlerJones,Markl}.
in~\cite{EnKoszulDuality}, we establish that this Koszul duality result at the homology level
has a counterpart at the chain level. Hence, we speak about a Koszul duality result
when we deal with $E_n$-operads.

The first reference addressing applications of homotopical algebra to algebras over operads
in the setting of unbounded dg-modules over a ring
is the article~\cite{HinichHomotopy}.
The definition of a homology theory with trivial coefficients, in terms of Quillen's homotopical algebra, was given earlier in~\cite{GetzlerJones}
in the context of non-negatively graded dg-modules over a field of characteristic $0$.
The authors of~\cite{GetzlerJones} also define an explicit cofibrant replacement functor on the category of $\POp$-algebras,
for any operad $\POp$,
by using the (already alluded to) operadic version of the cobar construction $B^c(\DOp)$.

In~\cite{CylinderOperads} we explain how to extend the applications of operadic bar duality to the context of unbounded dg-modules over a ring.
More specifically, we prove that the functorial cofibrant replacement of algebras over an operad, defined~\cite{GetzlerJones},
works in that framework provided that we restrict ourself to $\Sigma_*$-cofibrant operads
and to algebras which are cofibrant as dg-modules.

In this section, we briefly review the definition of these cofibrant replacement functors for algebras over operads,
and we explain the applications of the cofibrant replacement construction
to the homology of algebras over operads. We also explain the application of our general construction to $E_n$-operads
by using the Koszul duality result of~\cite{EnKoszulDuality}.

In this section (and in the next one), we assume that $\EOp_n$ is the particular $E_n$-operad used in the proof of the Koszul duality result of~\cite{EnKoszulDuality},
and which arises from a filtration of the Barratt-Eccles operad~\cite{BarrattEccles}.
The components of this operad $\EOp_n(r)$ are finite dimensional dg-modules,
and hence, can be dualized without care.
The operad $\EOp_n$ is also connected, and as a consequence, we have a well-defined cooperad structure
on the collection of dg-modules $\EOp_n(r)^{\vee}$
dual to $\EOp_n(r)$.

To begin this section, we say a few words about the operadic cobar construction~$B^c(\DOp)$.

\subsubsection{On operadic cobar constructions}\label{KoszulDualityReview:CobarConstruction}
The operadic cobar construction is a functor which associates a connected operad $B^c(\DOp)$
to any cooperad $\DOp$.
We do not really need the explicit definition of $B^c(\DOp)$.
Just recall that $B^c(\DOp)$ is a cofibrant as an operad whenever the cooperad $\DOp$ is $\Sigma_*$-cofibrant,
and any morphism $\phi: B^c(\DOp)\rightarrow\POp$
towards an operad $\POp$
is fully determined by a homogeneous morphism $\theta: \DOp\rightarrow\POp$, of degree $-1$,
vanishing on $\DOp(1)$,
and satisfying the equation:
\begin{equation}
\delta(\theta)\left\{\vcenter{\xymatrix@H=6pt@W=4pt@M=2pt@R=8pt@C=4pt{ i_1\ar[dr] & \cdots & i_r\ar[dl] \\
& *+<6pt>[F]{c}\ar[d] & \\
& 0 & }}\right\}
= \sum_{\substack{\tau\in\Theta_2(I)\\ \rho_{\tau}(\gamma)}}
\pm\lambda_*\left\{\vcenter{\xymatrix@H=6pt@W=4pt@M=2pt@R=8pt@C=4pt{ & i_*\ar[dr] & \cdots & i_*\ar[dl] & \\
i_*\ar[drr] & \cdots & *+<6pt>[F]{\theta(c_*)}\ar[d] & \cdots & i_*\ar[dll] \\
&& *+<6pt>[F]{\theta(c_*)}\ar[d] && \\
&& 0 && }}\right\},
\end{equation}
for any $c\in\DOp(r)$,
where $\lambda_*$ refers to the composition operation of the operad $\POp$
and we adopt conventions similar to~\S\ref{Background:TwistedObjects} to represent the expansion
of the coproduct of a cooperad element $c\in\DOp(r)$.
We refer to~\cite[\S 3.7]{CylinderOperads}
for details on this matter.

We adopt the notation $\phi = \phi_{\theta}$ for the morphism associated to $\theta: \DOp\rightarrow\POp$
and we refer to $\theta$
as the twisting cochain associated to $\phi_{\theta}$.

For our $E_n$-operad, we have:

\begin{thm}[{see~\cite[Theorems A-B]{EnKoszulDuality}}]\label{KoszulDualityReview:KoszulDualityStatement}
The operad $\EOp_n$ has a cofibrant model of the form $\QOp_n = B^c(\DOp_n)$,
where $\DOp_n = \Lambda^{-n}\EOp_n^{\vee}$
is the $n$-fold operadic desuspension of the dual cooperad of $\EOp_n$
in dg-modules.
\end{thm}

\begin{proof}[Explanations]
The proof of this assertion in~\cite{EnKoszulDuality}
is divided into several steps.
The first step involves the definition of a morphism $\phi_n: B^c(\DOp_n)\rightarrow\COp$
towards the operad of commutative algebras $\COp$.
The weak-equivalence $\psi_n: B^c(\DOp_n)\xrightarrow{\sim}\EOp_n$
is defined in a second step
by a lifting process from the morphism~$\phi_n$. We refer to~\cite{EnKoszulDuality}
for details.
\end{proof}

\subsubsection{Cofibrant replacements arising from operadic cobar constructions}\label{KoszulDualityReview:CofibrantConstruction}
Suppose now we have a $\Sigma_*$-cofibrant connected operad $\POp$
together with a cofibrant model of the form $\QOp = B^c(\DOp)$,
where $\DOp$ is a $\Sigma_*$-cofibrant operad.
Let $\phi_{\theta}: B^c(\DOp)\xrightarrow{\sim}\POp$
be the augmentation
given with the definition of this cofibrant model.

To any $\POp$-algebra $A$, we associate a (connected) quasi-cofree $\DOp$-coalgebra $\Gamma_{\POp}(A) = (\DOp(A),\partial_{\alpha})$
and a quasi-free $\POp$-algebra of the form $R_A = (\POp(\Gamma_{\POp}(A)),\partial_{\gamma})$.
The homomorphism $\alpha$ which determines the twisting coderivation of the quasi-cofree $\DOp$-coalgebra $\Gamma_{\POp}(A) = (\DOp(A),\partial_{\alpha})$
maps an element $c(a_1,\dots,a_r)\in\DOp(A)$
to the evaluation of the operation $\theta(c)\in\POp(r)$
on $a_1\otimes\dots\otimes a_r\in A^{\otimes r}$.
The homomorphism $\gamma$ which determines the twisting derivation of the quasi-free $\POp$-algebra $R_A = (\POp(\Gamma_{\POp}(A)),\partial_{\gamma})$
is given by the composite
\begin{equation*}
\Gamma_{\POp}(A)\xrightarrow{\rho}\DOp(\Gamma_{\POp}(A))\xrightarrow{\theta(\Gamma_{\POp}(A))}\POp(\Gamma_{\POp}(A)),
\end{equation*}
where $\rho$ is yielded by the universal $\DOp$-coalgebra structure of the quasi-cofree $\DOp$-coalgebra $\Gamma_{\POp}(A)$.

Then:

\begin{prop}[{see~\cite[\S 2]{GetzlerJones} and~\cite[\S 4.2]{CylinderOperads}}]\label{KoszulDualityReview:CofibrantReplacements}
We have a weak-equivalence $\epsilon: R_A\xrightarrow{\sim} A$ induced on generators $\Gamma_{\POp}(A)\subset R_A$
by the augmentation of the cofree $\DOp$-coalgebra $\epsilon: \DOp(A)\rightarrow A$,
and this $\POp$-algebra $R_A$ defines a cofibrant replacement
of the $\POp$-algebra $A$.\qed
\end{prop}

We refer to~\cite[\S 2]{GetzlerJones} for the version of this statement in the context of $\NN$-graded dg-modules over a field characteristic $0$.
We refer to~\cite[\S 4.2]{CylinderOperads} for a generalized statement which works in the setting of (unbounded) dg-modules over a ring.
We just use the simplifying observation that any $\POp$-algebra $A$
is cofibrant as a dg-module
to get the claim of our proposition
from the result established in this second reference.

In~\cite[\S 1.3]{EnKoszulDuality}, we check that the cooperad $\DOp_n = \Lambda^{-n}\EOp_n^{\vee}$
associated to our $E_n$-operad $\EOp_n$
is $\Sigma_*$-cofibrant (like the operad $\EOp_n$ itself).
Hence the result of Proposition~\ref{KoszulDualityReview:CofibrantReplacements}
holds for algebras over our $E_n$-operad.

\subsubsection{The cofibrant replacement functor and endomorphism operads}\label{KoszulDualityReview:EndomorphismOperads}
Recall that the action of an operad $\POp$
on an algebra $A$
is determined by an operad morphism $\nabla: \POp\rightarrow\End_A$,
where $\End_A$, the endomorphism operad of $A$,
is the operad such that $\End_A(r) = \Hom_{\C}(A^{\otimes r},A)$.
The evaluation of an operation $p\in\POp(r)$
on a tensor $a_1\otimes\dots\otimes a_r\in A^{\otimes r}$
is defined by the evaluation of the map $\nabla(p)\in\Hom_{\C}(A^{\otimes r},A)$
on $a_1\otimes\dots\otimes a_r$.

Observe that the homomorphism $\alpha: \DOp(A)\rightarrow A$
in the definition of $\Gamma_{\POp}(A)$
is, by adjunction,
equivalent to the homomorphism $\alpha^{\sharp}: \DOp\rightarrow\End_A$
such that $\phi_{\alpha^{\sharp}}$
is the composite
\begin{equation*}
B^c(\DOp)\xrightarrow{\phi_{\theta}}\POp\xrightarrow{\nabla_A}\End_A
\end{equation*}
giving the action of $B^c(\DOp)$ on $A$ through $\phi_{\theta}$.
Hence the twisting derivation of~$\Gamma_{\POp}(A)$
is fully determined by the restriction of the $\POp$-algebra structure of~$A$
to~$B^c(\DOp)$,
and not by the $\POp$-algebra structure itself.

\subsubsection{Operadic homology}\label{KoszulDualityReview:OperadicHomology}
Let $\POp$ be an augmented operad.
The indecomposable quotient of a $\POp$-algebra $A$ is defined as the quotient of the dg-module $A$
under the image of the operations $p: A^{\otimes r}\rightarrow A$,
where $p$ ranges over the augmentation ideal of $\POp$.
Thus, we have:
\begin{equation*}
\Indec_{\POp} A = A/\im(\widetilde{\POp}(A)\rightarrow A).
\end{equation*}

For a $\Sigma_*$-cofibrant augmented operad $\POp$,
we have a homology theory $H^{\POp}_*(-)$
defined by the formula
\begin{equation*}
H^{\POp}_*(A) = H_*(\Indec_{\POp} R_A)
\end{equation*}
for any $\POp$-algebra $A$,
where $R_A$ is a cofibrant replacement of $A$.

In the case of a quasi-free $\POp$-algebra,
we have an identity $\Indec_{\POp}(\POp(C),\partial_{\gamma}) = C$.
Thus, when we apply the definition of the homology $H^{\POp}_*(A)$ to the functorial cofibrant replacement $R_A = (\POp(\Gamma_{\POp}(A)),\partial_{\alpha})$
of Proposition~\ref{KoszulDualityReview:CofibrantReplacements},
we obtain the following assertion:

\begin{prop}\label{KoszulDualityReview:KoszulHomology}
We have the relation
\begin{gather*}
H^{\POp}_*(A) = H_*(\Gamma_{\POp}(A)) = H_*(\DOp(A),\partial_{\alpha})
\intertext{for any $\POp$-algebra $A$, where we consider the quasi-cofree coalgebra}
\Gamma_{\POp}(A) = (\DOp(A),\partial_{\alpha})
\end{gather*}
associated to $A$
in~\S\ref{KoszulDualityReview:CofibrantConstruction}.
\qed
\end{prop}

This proposition is a generalization in the context of unbounded dg-modules
over a field of possibly positive characteristic of a result
of~\cite{GetzlerJones}.

In the special case $\POp = \EOp_n$, the result of this proposition returns:
\begin{equation*}
H^{\EOp_n}_*(A) = H_*(\DOp_n(A),\partial_{\alpha}) = H_*(\Lambda^{-n}\EOp_n^{\vee}(A),\partial_{\alpha}).
\end{equation*}
Before applying this construction to the cochain algebra of a sphere,
we go back to the general study of coalgebras $\DOp(C)$.

\section{Free complete algebras and duality}\label{FreeCompleteAlgebras}
Throughout the paper, we adopt the notation $\FF[d]$ for a dg-module of rank one
concentrated in (lower) degree $d$.
The reduced normalized cochain complex of a sphere $\bar{N}^*(S^d)$ is an instance
of dg-module of this form.
We more precisely have $\bar{N}^*(S^d) = \FF[-d]$ since we change the sign of the degree when we move from upper to lower graded dg-modules.

The purpose of this section is to study coderivations of cofree coalgebras $\DOp(C)$
over a dg-module of this form $C = \FF[m]$ with the aim to determine,
in the next section,
the complex $(\Lambda^{-n}\EOp_n^{\vee}(A),\partial_{\alpha})$
associated to $A = \bar{N}^*(S^{n-m})$.
To make our description more conceptual, we define a dual complex in the setting of generalized power series algebras
rather than in a coalgebra setting. Therefore we study the duality process first.

From now on, we assume that $x$ is a variable of degree~$-m$.

\subsubsection{The free complete algebra on one variable over an operad}\label{FreeCompleteAlgebras:Definition}
In our presentation, we consider the free $\POp$-algebra $\POp(x)$ on a variable $x$,
or equivalently, the free $\POp$-algebra $\POp(\FF x)$
on the free dg-module of rank one spanned by $x$.
By definition, we have an identity $\POp(x) = \bigoplus_{r=0}^{\infty} (\POp(r)\otimes\FF x^{\otimes r})_{\Sigma_r}$.
To define the free complete $\POp$-algebra $\widehat{\POp}(x)$,
we just replace the sum in the expansion of the free $\POp$-algebra by a product (formed within the category of dg-modules):
$\widehat{\POp}(x) = \prod_{r=0}^{\infty} (\POp(r)\otimes\FF x^{\otimes r})_{\Sigma_r}$.
The free complete $\POp$-algebra
is equipped with the topology defined by the nested sequence of dg-submodules
$\widehat{\POp}{}^{(>s)}(x) = \prod_{r=s+1}^{\infty} (\POp(r)\otimes\FF x^{\otimes r})_{\Sigma_r}$.

\subsubsection{Duality}\label{FreeCompleteAlgebras:Duality}
By~\cite[Proposition 1.2.18]{KoszulOperads},
the dual in dg-modules of a connected cooperad $\DOp$
always forms a connected operad.
Now,
if $\DOp$ is equipped with a free $\Sigma_*$-structure,
then the dual of the connected cofree coalgebra $\DOp(C)$ on $C = \FF[m]$
can be identified with a complete free algebra $\widehat{\POp}(x)$ over the operad $\POp = \DOp^{\vee}$,
where $x$ is a variable of degree $-m$,
because we have in this case
$((\DOp(r)\otimes\FF[m]^{\otimes r})_{\Sigma_r})^{\vee}
\simeq(\DOp(r)^{\vee}\otimes\FF[-m]^{\otimes r})^{\Sigma_r}
\simeq(\DOp(r)^{\vee}\otimes\FF[-m]^{\otimes r})_{\Sigma_r}$
and the duality transforms the sum into a product.
In the sequel,
we are rather interested in a converse duality operation,
from free complete algebras over operads
to coalgebras over cooperads.

The dual~$\DOp = \POp^{\vee}$ of a connected operad~$\POp$ inherits a cooperad structure when we assume that each dg-module $\POp(r)$
is degree-wise finite dimensional over the ground field $\FF$.

In the sequel, we apply the duality of dg-modules to the free complete $\POp$-algebra $\widehat{\POp}(x)$.
In this case, we use the notation $\widehat{\POp}(x)^{\vee}$
to refer to the continuous dual with respect to the topology of the module $\widehat{\POp}(x)$.
Thus, we set $\widehat{\POp}(x)^{\vee} = \colim_{s} \{\widehat{\POp}(x)/\widehat{\POp}{}^{(>s)}(x)\}^{\vee}$.
We have the following result:

\begin{prop}\label{FreeCompleteAlgebras:DualityIsomorphism}
Let $x$ denote a homogeneous variable of degree $-m$ (as stated in the introduction of this section).
If the operad $\POp$ is equipped with a free $\Sigma_*$-structure and each dg-module $\POp(r)$
is degree-wise finite dimensional over the ground field $\FF$,
then we have an isomorphism $\widehat{\POp}(x)^{\vee}\simeq\DOp(\FF[m])$,
where $\DOp = \POp^{\vee}$
is the dual cooperad of~$\POp$.\qed
\end{prop}

Let $\gamma^{\sharp}\in\widehat{\POp}(x)^{\vee}$ denote the continuous homomorphism
associated to an element $\gamma = c(x^{\vee},\dots,x^{\vee})\in\DOp(\FF[m])$, $c\in\DOp(r)$,
where we use the notation $x^{\vee}$
for the canonical generator of the $\FF$-module~$\FF[m]$.
This homomorphism can be defined termwise by the relation $\gamma^{\sharp}(p(x,\dots,x)) = 0$, for $p\in\POp(s)$, $s\not=r$,
and the duality formula
\begin{equation*}
\gamma^{\sharp}(p(x,\dots,x)) = \sum_{w\in\Sigma_r} c(w\cdot p),
\end{equation*}
for $p\in\POp(r)$.

The relation $\widehat{\POp}(x)^{\vee}\simeq\DOp(\FF[m])$ is converse to the relation $\DOp(\FF[m])^{\vee}\simeq\widehat{\POp}(x)$
of~\S\ref{FreeCompleteAlgebras:Duality}
for~$\DOp = \POp^{\vee}$.
By suspension, we deduce from these duality relations an isomorphism between the suspended dg-module $\Sigma^{-m}\widehat{\POp}(x)$
and the dg-module of homomorphisms $\alpha: \DOp(\FF[m])\rightarrow\FF[m]$.

Recall that we also have a one-to-one correspondence between homomorphisms~$\alpha: \DOp(\FF[m])\rightarrow\FF[m]$
and coderivations~$\partial_{\alpha}: \DOp(\FF[m])\rightarrow\DOp(\FF[m])$.
This bijection defines an isomorphism of dg-modules when we equip the collection of coderivations $\partial_{\alpha}: \DOp(\FF[m])\rightarrow\DOp(\FF[m])$
with the natural dg-module structure inherited from the internal-hom of the category of dg-modules.
This correspondence has the following counterpart on the free complete $\POp$-algebra $\widehat{\POp}(x)$:

\begin{prop}\label{FreeCompleteAlgebras:ContinuousDerivations}
For any free complete $\POp$-algebra $\widehat{\POp}(x)$
on one variable $x$ of degree $m$,
we have a one-to-one correspondence
between:
\begin{enumerate}
\item
the elements $\omega\in\widehat{\POp}(x)$,
\item
the derivations of $\POp$-algebras $\partial: \widehat{\POp}(x)\rightarrow\widehat{\POp}(x)$
which are continuous with respect to the topology of $\widehat{\POp}(x)$.
\end{enumerate}
The degree of a derivation $\partial: \widehat{\POp}(x)\rightarrow\widehat{\POp}(x)$ is equal to the degree of the associated element $\omega$ in $\widehat{\POp}(x)$,
and our correspondence defines an isomorphism of dg-modules when we equip the collection of continuous derivations
with the natural dg-module structure inherited from the internal-hom of the category of dg-modules.
\end{prop}

In the sequel, we adopt the notation $\partial_{\omega}: \widehat{\POp}(x)\rightarrow\widehat{\POp}(x)$
for the continuous derivation associated to an element $\omega\in\widehat{\POp}(x)$.

\begin{proof}
In one direction, we retrieve the element $\omega$ associated to a derivation $\partial_{\omega}$
by setting $\omega = \partial_{\omega}(x)$.
To check that this definition gives a one-to-one correspondence, we essentially observe that a derivation $\partial$
satisfies the relation
\begin{equation*}
\partial(p(x,\dots,x)) = \sum_{i=1}^{r} p(x,\dots,\omega,\dots,x),
\end{equation*}
by definition, where we set $\omega = \partial(x)$. The sum runs over all inputs of the operation~$p\in\POp(r)$.
The continuity assumption implies that the derivation $\partial$
is determined on the whole dg-module~$\widehat{\POp}(x)$
by this term-wise relation.
\end{proof}

From the explicit construction of the derivation $\partial_{\omega}: \widehat{\POp}(x)\rightarrow\widehat{\POp}(x)$
in the proof of this proposition,
we obtain:

\begin{prop}\label{FreeCompleteAlgebras:DerivationAdjunction}
Let $\omega\in\widehat{\POp}(x)$.
Let $\omega^{\sharp}: \DOp(\FF[m])\rightarrow\FF[m]$
be the homomorphism associated to $\omega$
by the duality isomorphism $\widehat{\POp}(x)\simeq\DOp(\FF[m])^{\vee}$
of~\S\ref{FreeCompleteAlgebras:Duality}.

The adjoint homomorphism of the derivation $\partial_{\omega}: \widehat{\POp}(x)\rightarrow\widehat{\POp}(x)$
corresponds under the duality isomorphism $\widehat{\POp}(x)^{\vee}\simeq\DOp(\FF[m])$
to the coderivation~$\partial_{\omega^{\sharp}}: \DOp(\FF[m])\rightarrow\DOp(\FF[m])$
associated to~$\omega^{\sharp}$.
\end{prop}

\begin{proof}
Formal from the expression of the derivation~$\partial_{\omega}$ in the proof of proposition~\ref{FreeCompleteAlgebras:ContinuousDerivations}
of the coderivation~$\partial_{\omega^{\sharp}}$ in~\S\ref{Background:TwistedObjects},
and from the expression of the duality isomorphisms.
\end{proof}

The construction of quasi-free algebras over an operad can be adapted to free complete algebras.
In that situation,
we consider twisted dg-modules of the form $(\widehat{\POp}(x),\partial_\omega)$
where the twisting morphism $\partial_\omega$ is a continuous derivation of $\widehat{\POp}(x)$
associated to some element $\omega\in\widehat{\POp}(x)$.
As in the context of quasi-free algebras,
the derivation relation implies that the equation of twisting homomorphisms $\delta(\partial_{\omega}) + \partial_{\omega}^2 = 0$
is fulfilled for $\partial_\omega$
if and only if it holds on the generating element $x$
of $\widehat{\POp}(x)$.
Because of the relation $\omega = \partial_{\omega}(x)$,
this equation amounts to the identity $\delta(\omega) + \partial_{\omega}(\omega) = 0$
in $\widehat{\POp}(x)$.

The purpose of the next paragraphs is to explain that the composition structure on free complete algebras, which we consider in the introduction of this article,
represents the composition operation $\partial_{\alpha}(\beta)$,
where $\alpha,\beta\in\widehat{\POp}(x)$,
and we check that this composition structure can be used to characterize twisting derivations on $\widehat{\POp}(x)$.

\subsubsection{The composition structure}\label{FreeCompleteAlgebras:CompositionStructure}
To begin with, we revisit the definition of this composition structure in the general case
of a free complete algebra $\widehat{\POp}(x)$
over an operad $\POp$.
We adopt a formalism of~\cite{ChapotonLivernetHopf} (we also refer to~\cite{KapranovManin} for a variant of the construction).
We first define the composite of elements of homogeneous order $p(x,\dots,x),\break q(x,\dots,x)\in\widehat{\POp}(x)$
by a sum of operadic substitution operations
\begin{equation*}
p(x,\dots,x)\circ q(x,\dots,x) = \sum_{i=1}^{r} p(x,\dots,q(x,\dots,x),\dots,x)
\end{equation*}
running over all entries of $p(x,\dots,x)$.
This termwise composition operation is continuous,
and hence, extends to the completion $\lim_{r,s}\{\widehat{\POp}(x)/\widehat{\POp}{}^{(>r)}(x)\otimes\widehat{\POp}(x)/\widehat{\POp}{}^{(>s)}(x)\}$.
For our purpose, we just consider the restriction of the obtained composition operation
to the dg-module $\widehat{\POp}(x)\otimes\widehat{\POp}(x)
= \{\lim_{r}\widehat{\POp}(x)/\widehat{\POp}{}^{(>r)}(x)\}\otimes\{\lim_{s}\widehat{\POp}(x)/\widehat{\POp}{}^{(>s)}(x)\}$,
in order to obtain the operation:
\begin{equation*}
\circ: \widehat{\POp}(x)\otimes\widehat{\POp}(x)\rightarrow\widehat{\POp}(x).
\end{equation*}
Note that this composition operation
decreases degrees by $m$,
and becomes degree preserving
only when we form $m$-desuspensions of the dg-module $\widehat{\POp}(x)$.

This composition operation $\circ$ is, after desuspension, an instance of a pre-Lie structure,
a product satisfying the identity:
\begin{equation*}
(a\circ b)\circ b = a\circ(b\circ b),
\end{equation*}
for any pair of variables $(a,b)$. The notion of pre-Lie algebra was introduced in~\cite{Gerstenhaber}
in deformation theory. We refer to~\cite{ChapotonLivernetPreLie} for further historical references
on pre-Lie algebras, and a study of pre-Lie algebras from an operadic viewpoint.
Note that the definition of this reference is formulated in a characteristic zero setting.
The pre-Lie relation is equivalent to the identity $(a\circ b)\circ c - (a\circ c)\circ b = a\circ(b\circ c) - a\circ(c\circ b)$
in this setting, but this equivalence fails when our ground ring is a field of characteristic $2$.

The next lemma gives our motivation to introduce the pre-Lie composition operation:

\begin{lemm}\label{FreeCompleteAlgebras:CompositeDerivations}
For any pair $\alpha,\beta\in\widehat{\POp}(x)$,
we have the identity $\partial_{\alpha}(\beta) = \alpha\circ\beta$
in~$\widehat{\POp}(x)$.
\end{lemm}

\begin{proof}
This lemma is obvious from the explicit expression of the derivation $\partial_{\omega}$
associated to any $\omega\in\widehat{\POp}(x)$
in the proof of proposition~\ref{FreeCompleteAlgebras:ContinuousDerivations}.
\end{proof}

This lemma implies:

\begin{prop}\label{FreeCompleteAlgebras:TwistingDerivations}
The derivation $\partial_{\omega}: \widehat{\POp}(x)\rightarrow\widehat{\POp}(x)$
associated to any $\omega\in\widehat{\POp}(x)$
defines a twisting derivation on $\widehat{\POp}(x)$
if and only if we have the identity $\delta(\omega) + \omega\circ\omega = 0$
in $\widehat{\POp}(x)$.\qed
\end{prop}

In the next sections, we use the identity of Lemma~\ref{FreeCompleteAlgebras:CompositeDerivations}
to obtain a compact expression of the derivation $\partial_{\omega}$
associated to any element $\omega\in\widehat{\POp}(x)$
in Proposition~\ref{FreeCompleteAlgebras:ContinuousDerivations}.

Before proceeding to the study of the chain complex $C^{\EOp_n}_*(\bar{N}^*(S^{n-m}))$,
we record some applications of the pre-Lie composition structure
for the definition of morphisms
on the operad of $L_\infty$-algebras.

\subsubsection{On the operad of $L_\infty$-algebras}\label{FreeCompleteAlgebras:LinfinityOperad}
First of all,
recall that this operad is defined by the cobar construction $\LOp_{\infty} = \Lambda^{-1} B^c(\COp^{\vee})$,
where $\COp^{\vee}$ is the dual cooperad of the operad of commutative algebras.
The operad $\LOp_{\infty}$
comes equipped with a weak-equivalence $\epsilon: \LOp_{\infty}\xrightarrow{\sim}\LOp$,
where $\LOp$ is the operad of Lie algebras.
Hence,
we have an identity $H_*(\LOp_{\infty}) = \LOp$
and the homology $H_*(\LOp_{\infty})$
is generated as an operad by an operation $\lambda\in\LOp(2)$
satisfying the identities of a Lie bracket.
This observation is used in~\S\ref{HomotopyInvariance}.
For the moment,
we just want to record the following general proposition:

\begin{prop}\label{FreeCompleteAlgebras:LinfinityMorphisms}
If the operad $\POp$ is equipped with a free $\Sigma_*$-structure,
then we have a bijection between:
\begin{enumerate}
\item
the operad morphisms $\phi: \Lambda\LOp_{\infty}\rightarrow\Lambda^m\POp$;
\item
and the elements $\omega\in\widehat{\POp}(x)$,
of degree $-1-m$, and such that the relation $\delta(\omega) + \omega\circ\omega = 0$
holds in $\widehat{\POp}(x)$,
where we still assume that $x$ is a variable of degree~$-m$.
\end{enumerate}
\end{prop}

\begin{proof}
The cobar construction commutes with operadic suspensions.
Therefore we have an identity $\Lambda\LOp_\infty = B^c(\COp^{\vee})$.

Since $\COp^{\vee}(r) = \FF$ for all $r>0$,
the homomorphism $\theta: \COp^{\vee}\rightarrow\Lambda^m\POp$
associated to $\phi = \phi_{\theta}$
is equivalent to a collection of invariant elements in $\Lambda^m\POp(r)^{\Sigma_r}$.
We use the assumption about the $\Sigma_*$-structure of $\POp$
and the existence of an isomorphism $\POp(r)^{\Sigma_r}\simeq\POp(r)_{\Sigma_r}$
to identify this collection with an element of the form $\omega\in\widehat{\POp}(x)$.
We easily check that the equation of~\S\ref{KoszulDualityReview:CobarConstruction}
for the twisting cochain $\theta$
amounts to the equation $\delta(\omega) + \omega\circ\omega = 0$
for the element $\omega$
and the proposition follows.
\end{proof}

\section{The Koszul complex on the cochain algebras of spheres}\label{SphereKoszulDuality}

The aim of this section is to make explicit the complex $(\DOp_n(A),\partial_{\alpha})$
computing the homology~$H^{\EOp_n}_*(A)$
for the cochain algebra of a sphere $A = \bar{N}^*(S^{n-m})$.

We still consider the particular $E_n$-operads $\EOp_n$
arising from a filtration
\begin{equation}
\EOp_1\subset\cdots\subset\EOp_{n-1}\subset\EOp_n\subset\cdots\subset\EOp_\infty = \EOp
\end{equation}
of the chain Barratt-Eccles operad $\EOp$.
We do not need to review the definition of these operads $\EOp_n$
and of the Barratt-Eccles operad $\EOp$.
We are simply going to explain a representation of the $E_n$-algebra structure of $\bar{N}^*(S^{n-m})$
in terms of structures attached to our $E_n$-operads.
This representation will give the form of the complex $(\DOp_n(A),\partial_{\alpha})$
associated to $A = \bar{N}^*(S^{n-m})$.

We review the definition of the $\EOp_\infty$-algebra structure of $\bar{N}^*(S^{n-m})$ first,
before explaining applications of the Koszul duality result of~\cite{EnKoszulDuality}.

\subsubsection{The cochain algebra of spheres}\label{SphereKoszulDuality:SphereCochains}
Recall that the (reduced normalized) cochain complex of a $d$-sphere $S^d$
is identified with the free $\FF$-module of rank $1$
concentrated in (lower) degree $-d$:
\begin{equation*}
\bar{N}^*(S^d) = \FF[-d].
\end{equation*}
The associated endomorphism operad $\End_{\bar{N}^*(S^d)}$
is identified with the $d$-fold operadic desuspension of the commutative operad $\COp$
because we have identities:
\begin{equation*}
\End_{\bar{N}^*(S^d)}(r) = \Hom_{\C}(\FF[-d]^{\otimes r},\FF[-d]) = \FF[rd-d] = \Lambda^{-d}\COp(r).
\end{equation*}
Hence,
the $E_\infty$-structure of $\bar{N}^*(S^d)$
is determined by an operad morphism $\nabla_d: \EOp_\infty\rightarrow\Lambda^{-d}\COp$.

The work \cite{BergerFresse}
gives an explicit representation of this morphism for the chain Barratt-Eccles operad.
For our purpose, we use that $\nabla_d$
is identified with a composite
\begin{equation*}
\EOp_\infty\xrightarrow{\sigma}\Lambda^{-1}\EOp_\infty\xrightarrow{\sigma}\cdots\xrightarrow{\sigma}\Lambda^{-m}\EOp_\infty\xrightarrow{\epsilon}\Lambda^{-m}\COp,
\end{equation*}
where $\epsilon$ denotes the augmentation of the Barratt-Eccles operad, while $\sigma$ denotes a new operad morphism $\sigma: \EOp_\infty\rightarrow\Lambda^{-1}\EOp_\infty$
discovered in \cite{BergerFresse} and defined by an explicit formula. (For simplicity, we omit the application of operadic suspension
functors to these morphisms in the expression of our composite.)

In~\cite[Observation 0.1.6]{EnKoszulDuality}, we observe that $\sigma: \EOp_\infty\rightarrow\Lambda^{-1}\EOp_\infty$
restricts to an operad morphism $\sigma: \EOp_n\rightarrow\Lambda^{-1}\EOp_{n-1}$,
for each $n\geq 2$,
and we prove:

\begin{thm}[{see~\cite[Theorems A-B]{EnKoszulDuality}}]\label{SphereKoszulDuality:SuspensionKoszulDual}
The weak-equivalences $\psi_n: B^c(\DOp_n)\xrightarrow{\sim}\EOp_n$,
already considered in Theorem~\ref{KoszulDualityReview:KoszulDualityStatement},
fit a commutative diagram
\begin{equation}\label{EnEmbeddingModel}
\vcenter{\xymatrix{ B^c(\DOp_{n-1})\ar[r]^{\sigma^*}\ar@{.>}[d]_{\psi_{n-1}}^{\sim} &
B^c(\DOp_n)\ar@{.>}[d]_{\psi_n}^{\sim} \\
\EOp_{n-1}\ar[r]_{\iota} & \EOp_n }},
\end{equation}
where $\iota: \EOp_{n-1}\hookrightarrow\EOp_n$
is the embedding morphism
and $\sigma^*$ is the morphism associated to $\sigma: \EOp_n\rightarrow\Lambda^{-1}\EOp_{n-1}$
by functoriality of the construction $B^c(\DOp_m) = B^c(\Lambda^{-m}\EOp_m^{\vee})$.
\end{thm}

We already briefly recalled (see Theorem~\ref{KoszulDualityReview:KoszulDualityStatement})
that the Koszul duality weak-equivalence $\psi_m: B^c(\DOp_m)\xrightarrow{\sim}\EOp_m$
associated to our $E_m$-operad $\EOp_m$
arises as a lifting
\begin{equation*}
\xymatrix{ & \EOp_m\ar@{^{(}->}[]!R+<4pt,0pt>;[r] & \EOp\ar[d]^{\epsilon} \\
B^c(\DOp_m)\ar@{.>}[ur]^{\psi_m}\ar[rr]_{\phi_m} && \COp }
\end{equation*}
of a certain morphism $\phi_m$.
The notation $\epsilon$ refers again to the augmentation of the Barratt-Eccles operad.

By cobar duality, we have an equivalence between the diagram (\ref{EnEmbeddingModel})
and the diagram:
\begin{equation}\label{EnSuspensionModel}
\vcenter{\xymatrix{ B^c(\EOp_n^{\vee})\ar[r]^{\iota^*}\ar@{.>}[d]_{\psi_n^{\sharp}}^{\sim} &
B^c(\EOp_{n-1}^{\vee})\ar@{.>}[d]_{\psi_{n-1}^{\sharp}}^{\sim} \\
\Lambda^n\EOp_n\ar[r]_{\sigma} & \Lambda^{n-1}\EOp_{n-1} }}.
\end{equation}
The notation $\iota^*$
now refers to the morphism associated to the embedding $\iota: \EOp_{n-1}\hookrightarrow\EOp_n$
by functoriality of the construction $B^c((-)^{\vee})$.
From this result,
we deduce:

\begin{thm}[{see~\cite[Theorem C]{EnKoszulDuality}}]\label{SphereKoszulDuality:SphereAlgebraModel}
For every $n>m$,
we have a commutative diagram
\begin{equation*}
\xymatrix@C=1em{ B^c(\Lambda^{-n}\EOp_{n}^{\vee})\ar[d]_-{\iota^*}\ar[rr]^-{\sim} &&
\EOp_n\ar[d]^-{\sigma}\ar@{^{(}->}[]!R+<4pt,0pt>;[r] & \EOp\ar[d]^-{\sigma} \\
\vdots\ar[d]_-{\iota^*} && \vdots\ar[d]^-{\sigma} & \vdots\ar[d]^-{\sigma} \\
B^c(\Lambda^{-n}\EOp_{m}^{\vee})\ar[r]_-{\simeq} &
\Lambda^{m-n}B^c(\Lambda^{-m}\EOp_{m}^{\vee})\ar[r]_-{\sim}\ar@/_1em/[drr]_{\Lambda^{m-n}\phi'_{m}} &
\Lambda^{m-n}\EOp_{m}\ar@{^{(}->}[]!R+<4pt,0pt>;[r] &
\Lambda^{m-n}\EOp\ar[d]^-{\sim} \\
&&& \Lambda^{m-n}\COp = \End_{\bar{N}^*(S^{n-m})} }
\end{equation*}
giving the action of the $E_n$-operad $B^c(\Lambda^{-n}\EOp_n^{\vee})$
on $\bar{N}^*(S^{n-m})$,
where $\phi'_m$ is a morphism homotopic to $\phi_m$ in the model category of operads.
\end{thm}

Note: one can easily observe that the cobar construction commutes with operadic suspensions.

To determine the form of the complex $C_{\EOp_n}^*(\bar{N}^*(S^{n-m}))$,
we essentially use the observation of~\S\ref{KoszulDualityReview:EndomorphismOperads}
and the representation of the morphism $B^c(\DOp_n)\rightarrow\End_{\bar{N}^*(S^{n-m})}$
supplied by this theorem.

By~\cite[Theorem 5.2.2]{CylinderOperads} or~\cite[Theorem C]{PropHomotopy},
any $\QOp$-algebras $(A,\phi^0)$ and $(A,\phi^1)$
with the same underlying dg-module $A$
but a different $\QOp$-structure determined by morphisms $\phi^0,\phi^1: \QOp\rightarrow\End_A$
are connected by a chain of weak-equivalences
of $\QOp$-algebras
when the morphisms $\phi^0,\phi^1$ are homotopic in the category of operads.
Hence,
in our study of the cochain algebra $\bar{N}^*(S^{n-m})$,
any choice of morphism $B^c(\DOp_n)\rightarrow\End_{\bar{N}^*(S^{n-m})}$
in the homotopy class of the composite
\begin{equation*}
B^c(\DOp_n) = B^c(\Lambda^{-n}\EOp_{n}^{\vee})
\xrightarrow{\iota^*} B^c(\Lambda^{-n}\EOp_{m}^{\vee})
\simeq\Lambda^{m-n} B^c(\Lambda^{-m}\EOp_{m}^{\vee})
\xrightarrow{\phi'_m}\Lambda^{m-n}\COp
\end{equation*}
will give the right result.
For that reason,
we can replace the morphism $\phi'_m$ in Theorem~\ref{SphereKoszulDuality:SphereAlgebraModel}
by $\phi_m$,
or any homotopy equivalent morphism $\phi''_m$.
By the way, the construction of~\cite[\S 1]{EnKoszulDuality} actually gives a morphism $\phi_m$
which is only well-characterized up to homotopy.

Let us review the definition of these morphisms $\phi_{m}$
before going further in our study.

\subsubsection{The augmentation on the cobar constructions}\label{SphereKoszulDuality:CobarAugmentations}
Each operad morphism $\phi_m: B^c(\DOp_m)\rightarrow\COp$
is, according to the recollections of~\S\ref{KoszulDualityReview:CobarConstruction},
determined by a homomorphism of $\Sigma_*$-modules
\begin{equation*}
\DOp_m = \Lambda^{-m}\EOp_{m}^{\vee}\xrightarrow{\theta_m}\COp
\end{equation*}
satisfying a certain equation.
In~\cite[\S 1.1]{EnKoszulDuality}, we observe that the definition of such a map $\theta_m$
amounts to the definition of a collection of elements $\omega_m(r)\in\EOp_m(r)_{\Sigma_r}$
of degree $m(r-1)-1$,
or equivalently to an element $\omega_m$
of degree $-1-m$
in the completed $\EOp_m$-algebra $\widehat{\EOp}_m(x) = \prod_{r=1}^{\infty} (\EOp_m(r)\otimes\FF x^{\otimes r})_{\Sigma_r}$,
where $x$ is again a variable of degree $-m$.
The element $\omega_m(r)$
simply represents the term of order~$r$ of~$\omega_m\in\widehat{\EOp}_m(x)$.

Recall that $\COp(r) = \FF$ for all $r>0$.
The equivalence between $\omega_m(r)$
and $\theta_m$
is given by the adjunction relation $\theta_m(\xi) = \sum_{w\in\Sigma_r}\xi(w\cdot\underline{\omega}{}_m(r))$,
for any $\xi\in\Lambda^{-m}\EOp_{m}^{\vee}(r)$,
where $\underline{\omega}{}_m(r)$
is any representative of $\omega_m(r)$
in $\EOp_m(r)$.
In the relation,
we use the isomorphism $\EOp_m(r)_{\Sigma_r}\xrightarrow{\simeq}\EOp_m(r)^{\Sigma_r}$
defined by the norm of the $\Sigma_r$-action (recall that $\Sigma_r$ acts freely on $\EOp(r)$).

By inspection of~\cite[Proposition 1.1.A]{EnKoszulDuality},
we immediately see that the equation of the homomorphism $\theta_m$
amounts to the equation
\begin{equation}
\delta(\omega_m) + \omega_m\circ\omega_m = 0
\end{equation}
in term of the element $\omega_m\in\widehat{\EOp}_m(x)$,
where $\circ: \widehat{\EOp}_m(x)\otimes\widehat{\EOp}_m(x)\rightarrow\widehat{\EOp}_m(x)$
refers to the composition operation of~\S\ref{FreeCompleteAlgebras:CompositionStructure}.

Hence, we have by Proposition~\ref{FreeCompleteAlgebras:TwistingDerivations}:

\begin{lemm}\label{SphereKoszulDuality:TwistingDerivation}
The element $\omega_m\in\widehat{\EOp}_m(x)$
associated to the homomorphism $\theta_m$
defining $\phi_m$
determines, for each $n\geq m$, a twisting derivation $\partial_m: \widehat{\EOp}_n(x)\rightarrow\widehat{\EOp}_n(x)$
such that $\partial_m(\xi) = \xi\circ\omega_m$, for any $\xi\in\widehat{\EOp}_n(x)$.\qed
\end{lemm}

Then we observe:

\begin{lemm}\label{SphereKoszulDuality:DerivationAdjunctionRelation}
The adjoint homomorphism of the derivation $\partial_m: \widehat{\EOp}_n(x)\rightarrow\widehat{\EOp}_n(x)$
corresponds, under the isomorphisms
\begin{equation*}
\Sigma^{-n}\widehat{\EOp}_n(x)^{\vee}\simeq\Sigma^{-n}\EOp_n^{\vee}(\FF[m])\simeq\Lambda^{-n}\EOp_n^{\vee}(\FF[m-n]) = \DOp_n(N^*(S^{n-m}))
\end{equation*}
to the twisting coderivation~$\partial_{\alpha}: \DOp_n(N^*(S^{n-m}))\rightarrow\DOp_n(N^*(S^{n-m}))$
of the complex~$C^{\EOp_n}_*(N^*(S^{n-m}))$
computing $H^{\EOp_n}_*(N^*(S^{n-m}))$.
\end{lemm}

\begin{proof}
In~\S\ref{KoszulDualityReview:EndomorphismOperads},
we observe that the homomorphism $\alpha: \DOp_n(\FF[m-n])\rightarrow\FF[m-n]$
defining the coderivation of $C^*_{\EOp_n}(N^*(S^{n-m}))$.
corresponds, by adjunction, to the twisting cochain $\DOp_n\rightarrow\End_{\FF[m-n]} = \End_{N^*(S^{n-m})}$
determining the action of the operad $B^c(\DOp_n)$
on $N^*(S^{n-m})$.
By~Theorem~\ref{SphereKoszulDuality:SphereAlgebraModel},
this twisting cochain is given by a composite:
\begin{equation*}
\DOp_n\xrightarrow{\iota^*}\DOp_m\xrightarrow{\theta_m}\Lambda^{m-n}\COp = \End_{\FF[m-n]},
\end{equation*}
where $\iota^*$ is, up to operadic suspension,
the dual of the embedding morphism $\iota: \EOp_m\rightarrow\EOp_n$.
If we apply the formula of $\theta_m$
in terms of the element $\omega_m\in\widehat{\EOp}_m(x)$,
then we obtain the identity
\begin{equation*}
\alpha(c(x^{\vee},\dots,x^{\vee})) = \sum_{w\in\Sigma_r} c(w\cdot\omega_m(r))
\end{equation*}
for any $c\in\DOp_n(r) = \EOp_n(r)^{\vee}$,
where we go back to the notation used around Proposition~\ref{FreeCompleteAlgebras:DualityIsomorphism}
for the elements of~$\DOp_n(\FF[m-n])$.
Thus
we immediately see that $\alpha$ agrees with the homomorphism associated to $\omega_m\in\widehat{\EOp}_n(x)$
in~\S\S\ref{FreeCompleteAlgebras:DualityIsomorphism}-\ref{FreeCompleteAlgebras:ContinuousDerivations},
and the adjunction relation between the coderivation $\partial_{\alpha}$
and the derivation $\partial_m(\xi) = \xi\circ\omega_m$ associated to $\omega_m$
follows from the assertion of Proposition~\ref{FreeCompleteAlgebras:DerivationAdjunction}.
\end{proof}

From this statement, we conclude:

\begin{thm}\label{SphereKoszulDuality:HomologyComplex}
We have an identity $H_*^{\EOp_n}(\bar{N}^*(S^{n-m})) = \Sigma^{-n} H_*(\widehat{\EOp}_n(x)^{\vee},\partial_m^{\vee})$,
where we consider the twisting derivation $\partial_m(\xi) = \xi\circ\omega_m$
of Lemma~\ref{FreeCompleteAlgebras:ContinuousDerivations}
and the continuous dual of the twisted complete $\EOp_m$-algebra $(\widehat{\EOp}_n(x),\partial_m)$.\qed
\end{thm}

The theorem of the introduction
follows from this result
and from:

\begin{thm}[{see~\cite{EinfinityBar,IteratedBar}}]\label{SphereKoszulDuality:TopologicalInterpretation}
We take a finite primary field $\FF = \FF_p$ of characteristic $p>0$ as ground ring.
For any space $X$, whose homology $H_*(X)$ is a degreewise finite dimensional,
we have an identity
\begin{equation*}
H_*^{\EOp_n}(\bar{N}^*(X)) = \Sigma^{-n} H^*(\widehat{\Omega^n X}),
\end{equation*}
where $\widehat{\Omega^n X}$
refers to the $p$-profinite completion of the $n$-fold loop space $\Omega^n X$ (of which we take the continuous cohomology).\qed
\end{thm}

We have already explained in this section that the result of Theorem~\ref{SphereKoszulDuality:HomologyComplex}
does not depend on the choice of the morphism~$\phi_m$
in a certain homotopy class,
because any such choice gives rise to weakly-equivalent models of the $E_n$-algebra underlying $N^*(S^{n-m})$.
We prove in the next section that a more general homotopy invariance property holds for Theorem~\ref{SphereKoszulDuality:HomologyComplex}.

\subsubsection{Remark: the cohomological analogue of Theorem~\ref{SphereKoszulDuality:HomologyComplex}}
We have a result similar to Lemma~\ref{SphereKoszulDuality:DerivationAdjunctionRelation}
for the complex $C^*_{\EOp_n}(A,A)$
computing the natural cohomology theory $H^*_{\EOp_n}(A,A)$
associated to $E_n$-operads.
As a byproduct, we also have a cohomological analogue of the result of Theorem~\ref{SphereKoszulDuality:HomologyComplex}.
At the complex level,
we obtain an identity of the form
\begin{equation*}
C^*_{\EOp_n}(\bar{N}^*(S^{n-m}),\bar{N}^*(S^{n-m})) = \Sigma^{m}(\widehat{\EOp}_m(x),\ad_m),
\end{equation*}
where we equip the dg-module $\widehat{\EOp}_m(x)$
with the twisting cochain such that $\ad_m(\xi) = \xi\circ\omega_m - \pm\omega_m\circ\xi$.

\section{Homotopical invariance}\label{HomotopyInvariance}
In the previous section, we have used the Barratt-Eccles operad to get a specific model of $E_n$-operad,
but we mention in the introduction of this article
that the result of our theorems holds for any choice of $E_m$-operad $\EOp_m$
and for any choice of morphism $\phi_m: B^c(\DOp_m)\rightarrow\COp$.
The goal of this section is to prove this generalized homotopy invariance property.

Recall that $H_*(\EOp_m)$ is identified with the operad of associative algebras $\AOp$ for $m=1$,
with the $m$-Gerstenhaber operad $\GOp_m$ for $m>1$.
The associative operad is generated as an operad by an associative product operation $\mu\in\AOp(2)$,
of degree $0$.
The $m$-Gerstenhaber operad is generated by an associative and commutative product operation $\mu\in\GOp_m(2)$
of degree $0$
and by a Lie operation $\lambda\in\GOp_m(2)$
of degree $m-1$
satisfying a graded version of the Poisson distribution relation.

For the moment, we still consider the $E_m$-operad $\EOp_m$ of~\cite{EnKoszulDuality}.
The morphisms $\phi_m: B^c(\DOp_m)\rightarrow\COp$
are characterized by the relation $\phi_m(\mu) = \mu$ and $\phi_m(\lambda) = 0$
at the homology level,
because we have the following result:

\begin{lemm}[{see~\cite[Theorem B]{OperadMaps}}]\label{HomotopyInvariance:HomotopicalMorphisms}
Any pair of morphisms $\phi_m^{\epsilon}: B^c(\DOp_m)\rightarrow\COp$, $\epsilon = 0,1$,
satisfying the relations $\phi_m^{\epsilon}(\mu) = \mu$ and $\phi_m^{\epsilon}(\lambda) = 0$ in homology
are left homotopic in the model category of operads.\qed
\end{lemm}

Note that the relation $\phi_m(\lambda) = 0$ is obvious for $m>1$ since $\COp$
is concentrated in degree $0$.
For the element $\omega_m\in\widehat{\EOp}_m(x)$
associated to $\phi_m$,
the relation $\phi_m(\mu) = \mu$
amounts to the identity $[\omega_m(2)] = \mu$
in $H_*(\EOp_m(2)_{\Sigma_2})$,
where $\omega_m(2)$ refers to the component of~$\omega_m$
of order $2$.
Recall that $\omega_m(1) = 0$ and the relation $\delta(\omega_m) + \omega_m\circ\omega_m = 0$
implies at order $2$
that $\omega_m(2)$ defines a cycle in $\EOp_m(2)_{\Sigma_2}$ (see~\cite[Proposition 1.1.A]{EnKoszulDuality}).

Thus the lemma has an obvious interpretation in term of the element $\omega_m$.

\medskip
Now we plan to extend the result of Theorem~\ref{SphereKoszulDuality:HomologyComplex}
to $E_n$-operads $\Xi_n$
which are $\Sigma_*$-cofibrant but not necessarily finitely generated in all degree
and for each arity.
In that context,
we consider the operad $\LOp_{\infty}$,
associated to the usual category of $L_{\infty}$-algebras,
defined by the cobar construction $\LOp_{\infty} = \Lambda^{-1} B^c(\COp^{\vee})$,
where $\COp^{\vee}$ is the dual cooperad of the operad of commutative algebras.
Recall that $H_*(\LOp_{\infty})$ is isomorphic to the operad of Lie algebras $\LOp$
which is generated, as an operad, by an operation $\lambda\in\LOp(2)$
satisfying the identities of a Lie bracket (see~\S\ref{FreeCompleteAlgebras:LinfinityOperad}).

In the next statement,
we use the operadic suspension of this operad $\Lambda\LOp_{\infty} = B^c(\COp^{\vee})$
for which the homology is generated, as an operad, by a Lie bracket operation $\lambda\in\LOp(2)$
of degree $2$.
Since we can not form the dual of $\Xi_m$
when $\Xi_m$ is not free of finite rank in all degree
and for each arity,
we can not use the relationship between twisting elements $\omega_m\in\widehat{\Xi}_m(x)$
and morphisms $\phi_m: B^c(\Lambda^{-m}\Xi_m^{\vee})\rightarrow\COp$,
but we still have:

\begin{prop}\label{HomotopyInvariance:MorphismRepresentation}
The definition of a morphism $\phi_m^{\sharp}: \Lambda\LOp_{\infty}\rightarrow\Lambda^m\Xi_m$,
amounts to the definition of an element $\xi_m\in\widehat{\Xi}_m(x)$,
of degree $-1-m$,
vanishing at order~$1$,
and satisfying $\delta(\xi_m) + \xi_m\circ\xi_m = 0$
in $\widehat{\Xi}_m(x)$.
Moreover,
we have the identity $\phi_m^{\sharp}(\lambda) = \lambda$ in homology
if and only if $[\xi_m(2)] = \mu$
in $H_*(\Xi_m(2)_{\Sigma_2})$,
where~$\xi_m(2)$ represents the term of order~$2$ of~$\xi_m$.
\end{prop}

\begin{proof}
The first part of the proposition is established in Proposition~\ref{FreeCompleteAlgebras:LinfinityMorphisms}.
The second assertion of the proposition is an obvious consequence of the relation,
established in the proof of Proposition~\ref{FreeCompleteAlgebras:LinfinityMorphisms},
between the homomorphism $\theta_m^{\sharp}$ which determines $\phi_m^{\sharp}$
and the element $\xi_m$.
\end{proof}

The morphisms $\phi_m^{\sharp}: \Lambda\LOp_{\infty}\rightarrow\Lambda^m\Xi_m$
corresponds by bar duality to the morphisms $\phi_m: B^c(\Lambda^{-m}\Xi_m^{\vee})\rightarrow\COp$
of~\S\ref{SphereKoszulDuality:CobarAugmentations}
whenever this correspondence makes sense (see~\cite{OperadMaps}).

\medskip
Suppose now we have an $E_m$-operad $\Xi_m$ together with a morphism
of the form of Proposition~\ref{HomotopyInvariance:MorphismRepresentation}
and let $\xi_m$
be the associated element in the free complete $\Xi_m$-algebra $\widehat{\Xi}_m(x)$.
We aim at comparing the pair $(\Xi_m,\xi_m)$
to a pair $(\EOp_m,\omega_m)$
formed from the Barratt-Eccles operad.
We already know that the operads $\EOp_m$ and $\EOp_m$
are connected by a chain of weak-equivalences:
\begin{equation*}
\EOp_m\xleftarrow[\sim]{\phi}\Pi_m\xrightarrow[\sim]{\psi}\Xi_m.
\end{equation*}
We can moreover assume that the morphism $\phi$ and $\psi$
are fibrations.
We have then:

\begin{lemm}\label{HomotopyInvariance:MorphismCorrespondence}
We have elements $\pi_m\in\widehat{\Pi}_m(x)$ and $\omega_m\in\widehat{\EOp}_m(x)$
such that $\phi(\pi_m) = \omega_m$ and $\psi(\pi_m) = \xi_m$
and satisfying the identity $[\pi_m(2)] = [\omega_m(2)] = \mu$
in homology,
where as usual $\pi_m(2)$ (respectively, $\omega_m(2)$) denotes the term of order $2$ of $\pi_m$ (respectively, $\omega_m$).
\end{lemm}

\begin{proof}
The idea is to define the components $\pi_m(r)$ of the element $\pi_m\in\widehat{\Pi}_m(x)$
by induction on $r$.
For this, we use that the equation $\delta(\pi_m) + \pi_m\circ\pi_m = 0$
amounts at order $r$
to an equation of the form:
\begin{equation}\label{ObstructionEquation}
\delta(\pi_m(r)) = \sum_{s+t=r-1}\bigl\{\sum_{i=1}^{r}\pi_m(s)\circ_i\pi_m(t)\bigr\},
\end{equation}
where we use the standard partial composite notation of the operadic substitution operation
occurring in the definition of the composition product $\circ: \widehat{\Pi}_m(x)\otimes\widehat{\Pi}_m(x)\rightarrow\widehat{\Pi}_m(x)$.
Thus,
whenever appropriate elements $\pi_m(s)$ such that $\psi(\pi_m(s)) = \xi_m(s)$
are defined for $s<r$,
the existence of $\pi_m(r)$
amounts to the vanishing of an obstruction in $H_*(\Pi_m(r)_{\Sigma_r})$.
Now
this obstruction does vanish just because $\psi$ induces an isomorphism in homology,
preserves operadic structures,
and carries (\ref{ObstructionEquation})
to the same equation for $\xi_m$.
Hence we are done with the definition of $\pi_m$.

To define $\omega_m$,
we just set $\omega_m = \phi(\pi_m)$.
\end{proof}

\begin{lemm}\label{HomotopyInvariance:ComplexEquivalences}
We have a chain of weak-equivalences
\begin{equation*}
(\widehat{\EOp}_m(x),\partial_m)^{\vee}\xleftarrow[\sim]{\phi}(\widehat{\Pi}_m(x),\partial_m)^{\vee}
\xrightarrow[\sim]{\psi}(\widehat{\Xi}_m(x),\partial_m)^{\vee}
\end{equation*}
between the continuous duals of the chain complexes
formed from the elements of Lemma~\ref{HomotopyInvariance:MorphismCorrespondence}.
\end{lemm}

\begin{proof}
The commutation of $\phi$ and $\psi$ with operad structures readily implies that $\phi$ and $\psi$
induce dg-module morphisms between the twisted complexes.
By the standard spectral sequence argument,
we obtain that these dg-module morphisms are weak-equivalences
since $\phi$ and $\psi$ induce weak-equivalences
at the level of cofree coalgebras.
\end{proof}

Since the complex $(\widehat{\EOp}_n(x),\partial_m)$, which we obtain from the Barratt-Eccles operad in Lemma~\ref{HomotopyInvariance:ComplexEquivalences},
fits the construction of~\S\ref{SphereKoszulDuality},
we conclude:

\begin{thm}\label{HomotopyInvariance:Statement}
Let $\Xi_1\subset\dots\subset\Xi_n\subset\dots$ be any nested collection of $\Sigma_*$-cofibrant operads
equivalent to the nested sequence of the chain operads of little cubes.
Suppose we have an operad morphism
\begin{gather*}
\phi_m^{\sharp}: \Lambda\LOp_\infty\rightarrow\Lambda^m\Xi_m
\intertext{satisfying $\phi_m^{\sharp}(\lambda) = \lambda$
at the homology level,
and consider the equivalent element}
\xi_m\in\widehat{\Xi}_m(x)\quad\text{such that}\quad\delta(\xi_m) + \xi_m\circ\xi_m = 0.
\end{gather*}
The conclusion of Theorem~\ref{SphereKoszulDuality:HomologyComplex} holds for the complex $(\widehat{\Xi}_n(x),\partial_m)$
formed from this element $\xi_m\in\widehat{\Xi}_m(x)$.\qed
\end{thm}

In this theorem, we do not have to assume that the operad $\Xi_m$ consists of finitely generated dg-modules.
In particular, we can apply the theorem to the $E_n$-operads considered in~\cite{KontsevichMotives}
formed by the semi-algebraic chain complex of the Fulton-MacPherson operad.

\end{document}